\documentclass[11pt,reqno]{amsart}
\usepackage{graphicx}
\usepackage{verbatim}
\usepackage{textcomp}
\usepackage{amssymb}
\usepackage{cite}
\usepackage{amsmath}
\usepackage{latexsym}
\usepackage{amscd}
\usepackage{amsthm}
\usepackage{mathrsfs}
\usepackage{xypic}
\usepackage{bm}
\usepackage{url}
\usepackage{hyperref}

\newtheorem{thm}{Theorem}[section]
\newtheorem{corr}[thm]{Corollary}
\newtheorem{lem}[thm]{Lemma}
\newtheorem{prop}[thm]{Proposition}

\theoremstyle{definition}

\newtheorem*{ack}{Acknowledgment}
\theoremstyle{remark}
\newtheorem{rem}{Remark}[section]
\numberwithin{equation}{section}
\begin{document}
\title[Gradient
estimates  for a nonlinear elliptic equation]{Gradient estimates
for a nonlinear elliptic equation on complete Riemannian manifolds}
\author{Bingqing Ma}
\address{Department of Mathematics, Henan Normal
University, Xinxiang 453007, P.R. China,}
\email{bqma@henannu.edu.cn }

\author{Guangyue Huang}
\address{Department of Mathematics, Henan Normal
University, Xinxiang 453007, P.R. China,}
\email{hgy@henannu.edu.cn }

\author{Yong Luo}
\address{School of mathematics and statistics, Wuhan university, Wuhan 430072, P.R. China,}
\email{yongluo@whu.edu.cn }

\thanks{The research of the authors is supported by NSFC Nos.
}

\maketitle

\begin{abstract}
In this short note, we consider gradient estimates for positive
solutions to the following nonlinear elliptic equation on a complete
Riemannian manifold:
$$\Delta u+cu^{\alpha}=0,$$
where $c, \alpha$ are two real constants and $c\neq 0$.

\end{abstract}

\section{Introduction}

It is well-known that for complete noncompact Riemannian manifolds with nonnegative Ricci curvature, Yau \cite{Yau75} has
proved that any positive or bounded solution to the equation
\begin{equation}\label{Int2}
\Delta u=0\end{equation}
 must be constant. In \cite{Brighton2013}, Brighton studied $f$-harmonic function on a mooth metric measure
space. That is, he consider positive solutions to the equation
\begin{equation}\label{Int3}
\Delta_f u=0\end{equation} and obtain some similar results to Yau's
under the Bakry-\'{E}mery Ricci curvature condition.

It is easy to see that the equation \eqref{Int2} can be seen as a special case of
\begin{equation}\label{Int1}
\Delta u+cu^{\alpha}=0
\end{equation}
with $c ,\alpha$ two real constants. In particular, if $c=0$ in \eqref{Int1}, then the equation \eqref{Int1} becomes \eqref{Int2}. If $c<0$ and $\alpha<0$, the equation \eqref{Int1} on a bounded smooth domain in $\mathbb{R}^n$ is known as the thin film equation, which describes a steady state of the thin film (see \cite{Guowei-2006}). For $c$ a function, the equation \eqref{Int1} is studied by Gidas and Spruck in \cite{GidasSpruck-1981} with $1\leq \alpha\leq\frac{n+2}{n-2}$ when $n>2$ and lather it is studied by Li in \cite{Li-1991} to achieve gradient estimates and Liouville type results with $1<\alpha<\frac{n}{n-2}$ when $n>3$. In particular, Li achieved a gradient estimate for positive solution of (\ref{Int1}) when $c$ is a positive constant and $1<\alpha<\frac{n}{n-2}$.

Therefore, it is natural to try to achieve gradient estimates for positive solutions to the nonlinear elliptic equation \eqref{Int1} with other $c\neq 0$ and $\alpha$. In this direction Yang in \cite{Yang2010} proved the following result:
\begin{thm}[Yang]\label{Yang's thm}
Let $M$ be a complete noncompact  Riemannian manifold of dimension $n$ without boundary. Let $B_p(2R)$ be a geodesic ball of radius $2R$ around $p\in M$. We denote $-K(2R)$ with $K(2R)\geq 0$ such that $Ric_{ij}(B_p(2R))\geq -Kg_{ij}$. Suppose that $u(x)$ is a positive smooth solution of the equation (\ref{Int1}) with $\alpha<0$. Then we have

(i) If $c>0$, then $u(x)$ satisfies the estimate
$$\frac{|\nabla u|^2}{u^2}+cu^{\alpha-1}\leq C(n,\alpha)\big(K+\frac{1}{R^2}(1+\sqrt{K}R\coth
(\sqrt{K} R))\big)$$
on $B_p(R)$ and $C(n,\alpha)$ is a positive constant which depends on $n, \alpha$.

(ii) If $c<0$, then $u(x)$ satisfies the estimate
$$\frac{|\nabla u|^2}{u^2}+cu^{\alpha-1}\leq C(n,\alpha)\big(|c|(\inf_{B_p(2R)}u)^{\alpha-1}+K+\frac{1}{R^2}(1+\sqrt{K}R\coth
(\sqrt{K} R))\big)$$
on $B_p(R)$ and $C(n,\alpha)$ is a positive constant which depends on $n, \alpha$.
\end{thm}
After studying Yang's argument carefully, we find in the case of $c>0$, the gradient estimate in (i) actually holds when $\alpha\leq 1$, that is we have
\begin{thm}\label{1-thm-0}
Let $M$ be a complete noncompact  Riemannian manifold of dimension $n$ without boundary. Let $B_p(2R)$ be a geodesic ball of radius $2R$ around $p\in M$. We denote $-K(2R)$ with $K(2R)\geq 0$ such that $Ric_{ij}(B_p(2R))\geq -Kg_{ij}$. Suppose that $u(x)$ is a positive smooth solution of the equation (\ref{Int1}) with $\alpha\leq 1$ and $c>0$. Then we have
$$\frac{|\nabla u|^2}{u^2}+cu^{\alpha-1}\leq C(n,\alpha)\big(K+\frac{1}{R^2}(1+\sqrt{K}R\coth
(\sqrt{K} R))\big)$$
on $B_p(R)$ and $C(n,\alpha)$ is a positive constant which depends on $n, \alpha$.
\end{thm}
The proof of the above theorem is the same as Yang's proof of theorem \ref{Yang's thm}, and we will only give a sketch of it in the appendix. As a corollary of the above theorem we have the following Liouville type result:

\begin{corr}
Let $M$ be a complete noncompact  Riemannian manifold of dimension $n$ without boundary. Suppose that the Ricci curvature of $M$ is nonnegative. Then there does not exist a positive solution to equation (\ref{Int1}) with $\alpha\leq 1$ and $c>0$.
\end{corr}
Suppose that $u(x)$ is a positive solution to equation (\ref{Int1}). Following Brighton's argument in \cite{Brighton2013} by choosing a test function $u^\epsilon (\epsilon\neq 0)$, we can also get the following gradient estimate to $u(x)$.
\begin{thm}\label{1-thm-1}
Let $(M, g)$ be an $n$-dimensional complete Riemannian manifold with
$R_{ij}(B_p(2R))\geq-Kg_{ij}$, where $K\geq0$ is a constant. If $u$ is a positive solution to \eqref{Int1} on $B_p(2R)$ with $c$ and
$\alpha$ satisfying one of the following two cases:

(1) $c<0$ and $\alpha>0$;

(2) $c>0$ and $\frac{n+2}{2(n-1)}<\alpha<\frac{2n^2+9n+6}{2n(n+2)}$ with $n\geq3$,\\
then we have for any $x\in B_p(R)$
\begin{equation}\label{1-thmformula-1}\aligned
|\nabla u(x)|\leq&C(n,\alpha)M\sqrt{K+\frac{1}{R^2}\left( 1+ \sqrt{K}R \coth
(\sqrt{K} R)\right)},
\endaligned\end{equation}
where $M=\sup\limits_{x\in B_p(2R)}u(x)$ and the positive constant $C(n,\alpha)$
depends only on $n,\alpha$.
\end{thm}
\begin{rem}
In case (2), compared with Li's gradient estimate in \cite{Li-1991}
our right range for $\alpha$ is bigger than $\frac{n}{n-2}$ when
$n\geq 13$.
\end{rem}

Letting $R\rightarrow \infty$ in \eqref{1-thmformula-1}, we obtain
the following gradient estimates on complete noncompact Riemannian
manifolds:

\begin{corr}\label{1-corr-1}
Let $(M^n,g)$ be an $n$-dimensional complete noncompact Riemannian
manifold with $R_{ij}\geq-Kg_{ij}$, where $K\geq0$ is a constant.
Suppose that $u$ is a positive solution to \eqref{Int1} such that
$c, \alpha$ satisfies one of the two cases given in Theorem \ref{1-thm-1}. Then
we have
\begin{equation}\label{1-corrformula-1}\aligned
|\nabla u|\leq&C(n,\alpha)M\sqrt{K},
\endaligned\end{equation}
where $M=\sup\limits_{x\in M}u(x)$.

\end{corr}

\begin{rem}\label{rem3}
Recently, using the ideas of Brighton in \cite{Brighton2013}, some Liouville type results are achieved to positive solutions of the nonlinear elliptic equation
$$\Delta u+au\log u=0$$
in \cite{HM-2017}(for more developments, see \cite{Qian-2017,HLZ-2018}), and for porous medium and fast
diffusion equations in \cite{HM-GeomDedicata2017}.

\end{rem}

\begin{ack}
Yong Luo would like to thank Dr. Linlin Sun for his stimulating discussions on this problem.
\end{ack}

\section{Proof of theorem \ref{1-thm-1}}

Let $h=u^\epsilon$, where $\epsilon\neq0$ is a constant to be
determined. Then we have
\begin{equation}\label{Proof1}\aligned
\Delta h=&\epsilon(\epsilon-1)u^{\epsilon-2}|\nabla
u|^2+\epsilon u^{\epsilon-1}\Delta u\\
=&\epsilon(\epsilon-1)u^{\epsilon-2}|\nabla u|^2-c\epsilon u^{\alpha+\epsilon-1}\\
=&\frac{\epsilon-1}{\epsilon} \frac{|\nabla h|^2}{h}-c\epsilon
h^{\frac{\alpha+\epsilon-1}{\epsilon}},
\endaligned\end{equation}
where in the second equality of \eqref{Proof1}, we used \eqref{Int1}.
Hence, we have
\begin{equation}\label{Proof2}\aligned
\nabla h\nabla\Delta h=&\nabla
h\nabla\Big(\frac{\epsilon-1}{\epsilon} \frac{|\nabla
h|^2}{h}-c\epsilon
h^{\frac{\alpha+\epsilon-1}{\epsilon}}\Big)\\
=&\frac{\epsilon-1}{\epsilon}\nabla h\nabla \frac{|\nabla
h|^2}{h}-c(\alpha+\epsilon-1)h^{\frac{\alpha+\epsilon-1}{\epsilon}}\frac{|\nabla
h|^2}{h}\\
=&\frac{\epsilon-1}{\epsilon h}\nabla h\nabla (|\nabla
h|^2)-\frac{\epsilon-1}{\epsilon}\frac{|\nabla
h|^4}{h^2}-c(\alpha+\epsilon-1)h^{\frac{\alpha+\epsilon-1}{\epsilon}}\frac{|\nabla
h|^2}{h}.
\endaligned\end{equation}
Applying \eqref{Proof1} and \eqref{Proof2} into the well-known
Bochner formula to $h$, we have
\begin{equation}\label{Proof3}\aligned
\frac{1}{2}\Delta|\nabla h|^2=&|\nabla^2h|^2+\nabla h\nabla\Delta
h+{\rm Ric} (\nabla h,\nabla
h)\\
\geq&\frac{1}{n}(\Delta h)^2+\nabla h\nabla\Delta
h-K|\nabla h|^2\\
=&\frac{1}{n}\Big(\frac{\epsilon-1}{\epsilon} \frac{|\nabla
h|^2}{h}-c\epsilon
h^{\frac{\alpha+\epsilon-1}{\epsilon}}\Big)^2+\frac{\epsilon-1}{\epsilon
}\frac{\nabla h}{h}\nabla
(|\nabla h|^2)\\
&-\frac{\epsilon-1}{\epsilon}\frac{|\nabla h|^4}{h^2}
-c(\alpha+\epsilon-1)h^{\frac{\alpha+\epsilon-1}{\epsilon}}\frac{|\nabla
h|^2}{h}-K|\nabla h|^2\\
=&\Big(\frac{(\epsilon-1)^2}{n\epsilon^2}-\frac{\epsilon-1}{\epsilon}\Big)\frac{|\nabla
h|^4}{h^2}-c\left[\frac{n+2}{n}(\epsilon-1)+\alpha\right]
h^{\frac{\alpha+\epsilon-1}{\epsilon}}\frac{|\nabla
h|^2}{h}\\
&+\frac{c^2\epsilon^2}{n}h^{\frac{2(\alpha+\epsilon-1)}{\epsilon}}+\frac{\epsilon-1}{\epsilon
}\frac{\nabla h}{h}\nabla (|\nabla h|^2)-K|\nabla h|^2.
\endaligned\end{equation}
By analyzing (\ref{Proof3}) we have the following lemmas.

\begin{lem}\label{22-Lemma1}
Let $u$ be a positive solution to \eqref{Int1} and
$R_{ij}\geq-Kg_{ij}$ for some nonnegative constant $K$. Denote
$h=u^\epsilon$ with $\epsilon\neq0$. If $c<0$ and $\alpha>0$, then there exists $\epsilon\in(0,1)$ such that
\begin{equation}\label{0-Lemma-1}\aligned
\frac{1}{2}\Delta|\nabla
h|^2\geq&\Big(\frac{(\epsilon-1)^2}{n\epsilon^2}-\frac{\epsilon-1}{\epsilon}
\Big)\frac{|\nabla
h|^4}{h^2}\\
&+\frac{\epsilon-1}{\epsilon }\frac{\nabla h}{h}\nabla (|\nabla
h|^2)-K|\nabla h|^2.
\endaligned\end{equation}
\end{lem}

\proof In (\ref{Proof3}), if $c<0$ and $\alpha>0$, we can choose $\epsilon \in (0,1)$ close enough to 1 such that $$-c\left[\frac{n+2}{n}(\epsilon-1)+\alpha\right]\geq0,$$
and then (\ref{0-Lemma-1}) follows directly.
\endproof

\begin{lem}\label{22-Lemma3}
Let $u$ be a positive solution to \eqref{Int1} and
$R_{ij}\geq-Kg_{ij}$ for some nonnegative constant $K$. Denote
$h=u^\epsilon$ with $\epsilon\neq0$. If $c>0$ and for a fixed $\alpha$, there
exist two positive constants $\epsilon,\delta$ such that
\begin{equation}\label{1-Lemma-1}\aligned
c\left[\frac{n+2}{n}(\epsilon-1)+\alpha\right]>0
\endaligned\end{equation}
and
\begin{equation}\label{1-Lemma-2}\aligned
\frac{c^2\epsilon^2}{n}-\frac{c}{\delta}\Big(\frac{n+2}{n}(\epsilon-1)+\alpha\Big)>0,
\endaligned\end{equation}
then we have
\begin{equation}\label{1-Lemma-3}\aligned
\frac{1}{2}\Delta|\nabla
h|^2\geq&\Big[\frac{(\epsilon-1)^2}{n\epsilon^2}-\frac{\epsilon-1}{\epsilon}
-\delta c\Big(\frac{n+2}{n}(\epsilon-1)+\alpha\Big)\Big]\frac{|\nabla
h|^4}{h^2}\\
&+\frac{\epsilon-1}{\epsilon }\frac{\nabla h}{h}\nabla (|\nabla
h|^2)-K|\nabla h|^2.
\endaligned\end{equation}
\end{lem}
\proof For a fixed point $p$, if there exists a positive constant
$\delta$ such that $h^{\frac{\alpha+\epsilon-1}{\epsilon}}\leq\delta\frac{|\nabla
h|^2}{h}$, according to \eqref{1-Lemma-1}, then \eqref{Proof3} becomes
\begin{equation}\label{Proof5}\aligned
\frac{1}{2}\Delta|\nabla
h|^2\geq&\Big[\frac{(\epsilon-1)^2}{n\epsilon^2}-\frac{\epsilon-1}{\epsilon}
-\delta c\Big(\frac{n+2}{n}(\epsilon-1)+\alpha\Big)\Big]\frac{|\nabla
h|^4}{h^2}\\
&+\frac{c^2\epsilon^2}{n}h^{\frac{2(\alpha+\epsilon-1)}{\epsilon}}+\frac{\epsilon-1}{\epsilon
}\frac{\nabla
h}{h}\nabla (|\nabla h|^2)-K|\nabla h|^2\\
\geq&\Big[\frac{(\epsilon-1)^2}{n\epsilon^2}-\frac{\epsilon-1}{\epsilon}
-\delta c\Big(\frac{n+2}{n}(\epsilon-1)+\alpha\Big)\Big]\frac{|\nabla
h|^4}{h^2}\\
&+\frac{\epsilon-1}{\epsilon }\frac{\nabla h}{h}\nabla (|\nabla
h|^2)-K|\nabla h|^2.
\endaligned\end{equation}
On the contrary, at the point $p$, if
$h^{\frac{\alpha+\epsilon-1}{\epsilon}}\geq\delta\frac{|\nabla h|^2}{h}$, then
\eqref{Proof3} becomes
\begin{equation}\label{Proof6}\aligned
\frac{1}{2}\Delta|\nabla
h|^2\geq&\Big(\frac{(\epsilon-1)^2}{n\epsilon^2}-\frac{\epsilon-1}{\epsilon}\Big)\frac{|\nabla
h|^4}{h^2}
+\Big[\frac{c^2\epsilon^2}{n}-\frac{c}{\delta}\Big(\frac{n+2}{n}(\epsilon-1)+\alpha\Big)\Big]
h^{\frac{2(\alpha+\epsilon-1)}{\epsilon}}\\
&+\frac{\epsilon-1}{\epsilon }\frac{\nabla h}{h}\nabla (|\nabla
h|^2)-K|\nabla h|^2\\
\geq&\Big\{\Big(\frac{(\epsilon-1)^2}{n\epsilon^2}-\frac{\epsilon-1}{\epsilon}\Big)
+\delta^2\Big[\frac{c^2\epsilon^2}{n}-\frac{c}{\delta}\Big(\frac{n+2}{n}(\epsilon-1)+\alpha\Big)\Big]\Big\}\frac{|\nabla
h|^4}{h^2}\\
&+\frac{\epsilon-1}{\epsilon }\frac{\nabla h}{h}\nabla (|\nabla
h|^2)-K|\nabla h|^2\\
\geq&\Big[\frac{(\epsilon-1)^2}{n\epsilon^2}-\frac{\epsilon-1}{\epsilon}
-\delta c\Big(\frac{n+2}{n}(\epsilon-1)+\alpha\Big)\Big]\frac{|\nabla
h|^4}{h^2}\\
&+\frac{\epsilon-1}{\epsilon }\frac{\nabla h}{h}\nabla (|\nabla
h|^2)-K|\nabla h|^2
\endaligned\end{equation} as long as
\begin{equation}\label{Proof7}
\frac{c^2\epsilon^2}{n}-\frac{c}{\delta}\Big(\frac{n+2}{n}(\epsilon-1)+\alpha\Big)>0.
\end{equation}

In both cases, \eqref{1-Lemma-3} holds always. We complete the proof of Lemma \ref{22-Lemma3}.\endproof

In order to obtain the upper bound of $|\nabla h|$ by using the maximum
principle, it is sufficient to choose the
coefficient of $\frac{|\nabla h|^4}{h^2}$ in \eqref{0-Lemma-1} and \eqref{1-Lemma-3} such that it is
positive. In case of Lemma \ref{22-Lemma3}, we need to choose appropriate $\epsilon, \delta$ such that
\begin{equation}\label{Proof8}
\frac{(\epsilon-1)^2}{n\epsilon^2}-\frac{\epsilon-1}{\epsilon}
-\delta c\Big(\frac{n+2}{n}(\epsilon-1)+\alpha\Big)>0.\end{equation}
Under the assumption of \eqref{1-Lemma-1}, the inequality
\eqref{1-Lemma-2} becomes
\begin{equation}\label{Proof10}
\delta>\frac{nc}{c^2
\epsilon^2}\Big(\frac{n+2}{n}(\epsilon-1)+\alpha\Big)
\end{equation}
and \eqref{Proof8} becomes \begin{equation}\label{Proof11}
\delta<\frac{\frac{(\epsilon-1)^2}{n\epsilon^2}
-\frac{\epsilon-1}{\epsilon}}{c\Big(\frac{n+2}{n}(\epsilon-1)+\alpha\Big)}.
\end{equation}
In order to ensure we can choose a positive $\delta$, from \eqref{Proof10} and  \eqref{Proof11}, we need choose an $\epsilon$ satisfying
\begin{equation}\label{add-Proof12}
\frac{nc}{c^2
\epsilon^2}\Big(\frac{n+2}{n}(\epsilon-1)+\alpha\Big)<\frac{\frac{(\epsilon-1)^2}{n\epsilon^2}
-\frac{\epsilon-1}{\epsilon}}{c\Big(\frac{n+2}{n}(\epsilon-1)+\alpha\Big)}.
\end{equation}
In particular, \eqref{add-Proof12} can be written as
\begin{equation}\label{add2-Proof12}\aligned
n^2\Big(\frac{n+2}{n}(\epsilon-1)+\alpha\Big)^2<&n\epsilon^2\Big(\frac{(\epsilon-1)^2}{n\epsilon^2}-\frac{\epsilon-1}{\epsilon}\Big)\\
=&(\epsilon-1)^2-n\epsilon(\epsilon-1),
\endaligned\end{equation}
which is equivalent to
\begin{equation}\label{Proof12}\aligned
{}[n^2+5n+3]\epsilon^2&+[2(\alpha-1)(n^2+2n)-(5n+6)]\epsilon\\
&+(\alpha-1)^2n^2-4(\alpha-1)n+3<0.
\endaligned\end{equation}
By a direct calculation, under the condition
\begin{equation}\label{Proof14}\aligned
\frac{-(n-4)-\sqrt{n^2+5n+3}}{2(n-1)}<\alpha-1<\frac{-(n-4)+ \sqrt{n^2+5n+3}}{2(n-1)},
\endaligned\end{equation}
we have
\begin{equation}\label{Proof13}\aligned
{}[2&(\alpha-1)(n^2+2n)-(5n+6)]^2-4[n^2+5n+3][(\alpha-1)^2n^2-4(\alpha-1)n+3]\\
=&4(\alpha-1)^2[(n^2+2n)^2-n^2(n^2+5n+3)]+4(\alpha-1)[4n(n^2+5n+3)\\
&-(n^2+2n)(5n+6)]+(5n+6)^2-12(n^2+5n+3)\\
=&4(\alpha-1)^2[-n^3+n^2]+4(\alpha-1)[-n^3+4n^2]+13n^2\\
=&n^2\Big\{-4(n-1)(\alpha-1)^2-4(n-4)(\alpha-1)+ 13\Big\}\\
>&0,
\endaligned\end{equation}
which shows the quadratic inequality \eqref{Proof12} with respect to $\epsilon$ has two real roots.

Now we are ready to prove the following proposition.
\begin{prop}\label{22-prop2}
Let $u$ be a positive solution to \eqref{Int1} and
$R_{ij}\geq-Kg_{ij}$ for some nonnegative constant $K$. If we choose
$c$ and $\alpha$ satisfies one of the following two cases:

(1) $c<0$ and $\alpha>0$;

(2) $c>0$ and
$\frac{n+2}{2(n-1)}<\alpha<\frac{2n^2+9n+6}{2n(n+2)}$ with $n\geq3$, Then we have
\begin{equation}\label{add-2-Lemma-2}\aligned
\frac{1}{2}\Delta|\nabla h|^2\geq&C_1(n,\alpha)\frac{|\nabla
h|^4}{h^2}-C_2(n,\alpha)\frac{\nabla h}{h}\nabla (|\nabla h|^2)-K|\nabla
h|^2,
\endaligned\end{equation}
where $C_1(n,\alpha)$ and $C_2(n,\alpha)$ are positive constants.

\end{prop}

\proof We prove this proposition case by case.

\noindent{\bf (i) The case of $c<0$ and $\alpha>0$.} In the proof of
Lemma \ref{22-Lemma1} we see that by choosing an
$\epsilon=\epsilon(n,\alpha)\in (0,1)$ such that
$\frac{n+2}{n}(\epsilon-1)+\alpha\geq 0$ we get the
\begin{equation}\aligned \frac{1}{2}\Delta|\nabla
h|^2\geq&\Big(\frac{(\epsilon-1)^2}{n\epsilon^2}-\frac{\epsilon-1}{\epsilon}
\Big)\frac{|\nabla
h|^4}{h^2}\\
&+\frac{\epsilon-1}{\epsilon }\frac{\nabla h}{h}\nabla (|\nabla
h|^2)-K|\nabla h|^2.
\endaligned\end{equation}
Then we see that $C_1(n,\alpha)=\frac{(\epsilon-1)^2}{n\epsilon^2}-\frac{\epsilon-1}{\epsilon}>0$ and $C_2(n,\alpha)=\frac{1-\epsilon}{\epsilon}>0$.

\noindent{\bf (ii) The case of  $c>0$ and $\frac{n+2}{2(n-1)}<\alpha<\frac{2n^2+9n+6}{2n(n+2)}$ when $n\geq 3$.} In this case,
\eqref{1-Lemma-1} is equivalent to
\begin{equation}\label{Proof16}
\epsilon>1-\frac{n\alpha}{n+2}.
\end{equation}
We can check
\begin{equation}\label{Proof18}
\frac{5n+6}{2(n^2+2n)}<\frac{-(n-4)+\sqrt{n^2+5n+3}}{2(n-1)}.
\end{equation}
Hence, when $n\geq3$, for any $\alpha$ satisfies
\begin{equation}\label{Proof17}
-\frac{n-4}{2(n-1)}<\alpha-1<\frac{5n+6}{2(n^2+2n)}
\end{equation}
which is equivalent to
\begin{equation}\label{add2-Proof17}
\frac{n+2}{2(n-1)}<\alpha<\frac{2n^2+9n+6}{2n(n+2)},
\end{equation}
then
\eqref{Proof16} is satisfied by choosing
 \begin{equation}\label{2-Lemma-3}
\epsilon:=\tilde{\epsilon}=\frac{(5n+6)-2(\alpha-1)(n^2+2n)}{2(n^2+5n+3)},
\end{equation}
and it is easy to check that $\epsilon \in (0,1)$.

In particular, we let
\begin{equation}\label{Proof20}
\delta=\tilde{\delta}:=\frac{1}{2}\Bigg[\frac{nc}{c^2
\tilde{\epsilon}^2}\Big(\frac{n+2}{n}(\tilde{\epsilon}-1)+\alpha\Big)+\frac{\frac{(\tilde{\epsilon}-1)^2}{n\tilde{\epsilon}^2}
-\frac{\tilde{\epsilon}-1}{\tilde{\epsilon}}}{c\Big(\frac{n+2}{n}(\tilde{\epsilon}-1)+\alpha\Big)}\Big],
\end{equation}
then \eqref{Proof7} and \eqref{Proof8} are satisfied and
\eqref{1-Lemma-3} becomes
\begin{equation}\label{2-Lemma-2}\aligned
\frac{1}{2}\Delta|\nabla h|^2\geq&\tilde{C_1}(n,\alpha)\frac{|\nabla
h|^4}{h^2}-\tilde{C_2}(n,\alpha)\frac{\nabla h}{h}\nabla (|\nabla h|^2)-K|\nabla
h|^2,
\endaligned\end{equation}
where positive constants $\tilde{C_1}(n,\alpha)$ and $\tilde{C_2}(n,\alpha)$ are given by
$$\tilde{C_1}(n,\alpha)=\frac{1}{2}\Big[\Big(\frac{(\tilde{\epsilon}-1)^2}{n\tilde{\epsilon}^2}
-\frac{\tilde{\epsilon}-1}{\tilde{\epsilon}}\Big)
-\frac{n}{\tilde{\epsilon}^2}\Big(\frac{n+2}{n}(\tilde{\epsilon}-1)+\alpha\Big)^2\Big],$$
$$\tilde{C_2}(n,\alpha)=\frac{4(\alpha-1)n(n+2)+n(2n+5)}{(5n+6)-4(\alpha-1)n(n+2)},$$
respectively.

We conclude the proof of Proposition \ref{22-prop2}.\endproof

Now, we are in a position to prove our Theorem \ref{1-thm-1}. Denote
by $B_p(R)$ the geodesic ball centered at $p$ with radius $R$. Let
$\phi$ be a cut-off function (see \cite{schoenyau}) satisfying
$\mathrm{supp}(\phi) \subset B_p(2R)$, $\phi |_{B_p(R)}=1$
and
\begin{equation}\label{2-Proof25}
\frac{|\nabla \phi|^2}{\phi}\leq \frac{C}{R^2},\end{equation}
\begin{equation}\label{2-Proof26}
-\Delta \phi \leq\frac{C}{R^2}\left(1+\sqrt{K}R\coth (\sqrt{K}
R)\right),
\end{equation}
where $C$ is a constant depending only on $n$. We define $G=\phi
|\nabla h|^2$ and will apply the maximum principle to $G$ on $B_p(2
R)$. Moreover, we assume $G$ attains its maximum at the point
$x_0\in B_p(2R)$ and assume $G(x_0)>0$ (otherwise the proof is
trivial). Then at the point $x_0$, it holds that
$$\Delta G\leq 0, \ \ \
\nabla (|\nabla h|^2)=-\frac{|\nabla h|^2}{\phi} \nabla \phi$$ and
\begin{equation}\label{2-Proof27}\aligned
0\geq& \Delta G\\
=&\phi\Delta(|\nabla h|^2)+|\nabla h|^2\Delta\phi+2\nabla\phi\nabla
|\nabla h|^2 \\
=&\phi\Delta(|\nabla h|^2)+\frac{\Delta\phi}{\phi}G
-2\frac{|\nabla\phi|^2}{\phi^2}G\\
\geq&2\phi\Big[C_1(n,\alpha)\frac{|\nabla h|^4}{h^2}-C_2(n,\alpha)\frac{\nabla
h}{h}\nabla (|\nabla h|^2)-K|\nabla h|^2\Big]\\
&+\frac{\Delta\phi}{\phi}G -2\frac{|\nabla\phi|^2}{\phi^2}G\\
=&2C_1(n,\alpha)\frac{G^2}{\phi
h^2}+2C_2(n,\alpha)\frac{G}{\phi}\nabla\phi\frac{\nabla
h}{h}-2KG+\frac{\Delta\phi}{\phi}G -2\frac{|\nabla\phi|^2}{\phi^2}G,
\endaligned\end{equation}
where, in the second inequality, the estimate \eqref{2-Lemma-2} is
used. Multiplying both sides of \eqref{2-Proof27} by $\frac{\phi}{G}$
yields
\begin{equation}\label{2-Proof28}\aligned
2C_1(n,\alpha)\frac{G}{ h^2}\leq-2C_2(n,\alpha)\nabla\phi\frac{\nabla
h}{h}+2\phi K- \Delta\phi+2\frac{|\nabla\phi|^2}{\phi}.
\endaligned\end{equation}
Using the Cauchy inequality
$$\aligned
-2C_2(n,\alpha)\nabla\phi\frac{\nabla
h}{h}\leq&2C_2(n,\alpha)|\nabla\phi|\frac{|\nabla
h|}{h}\\
\leq&\frac{C_2^2(n,\alpha)}{C_1(n,\alpha)}\frac{|\nabla\phi|^2}{\phi}+C_1(n,\alpha)
\frac{G}{h^2},
\endaligned$$
into \eqref{2-Proof28} yields
\begin{equation}\label{2-Proof29}\aligned
C_1(n,\alpha)\frac{G}{ h^2}\leq2\phi K-
\Delta\phi+\Big(2+\frac{C_2^2(n,\alpha)}{C_1(n,\alpha)}\Big)\frac{|\nabla\phi|^2}{\phi}.
\endaligned\end{equation}
Hence, for $x\in B_p(R)$, we have
\begin{equation}\label{2-Proof30}\aligned
C_1(n,\alpha)G(x)\leq&C_1(n,\alpha)G(x_0)\\
\leq&h^2(x_0)\left[2K+\frac{C(n,p)}{R^2}\left( 1+ \sqrt{K}R \coth
(\sqrt{K} R)\right)\right].
\endaligned\end{equation}
It shows that
\begin{equation}\label{2-Proof31}\aligned
|\nabla u|^2(x)\leq&C(n,\alpha)M^2\left[K+\frac{1}{R^2}\left( 1+
\sqrt{K}R \coth (\sqrt{K} R)\right)\right]
\endaligned\end{equation}
and hence,
\begin{equation}\label{2-Proof32}\aligned
|\nabla u(x)|\leq&C(n,\alpha)M\sqrt{K+\frac{1}{R^2}\left( 1+ \sqrt{K}R \coth
(\sqrt{K} R)\right)}.
\endaligned\end{equation}

We complete the proof of Theorem \ref{1-thm-1}.

\section{Appendix}
Here we give a sketch of the proof of theorem \ref{1-thm-0}. The interested readers can consult Yang's paper \cite{Yang2010} for details. Assume that $u(x)$ is a positive solution to (\ref{Int1}) with $c>0$ and $\alpha\leq 1$. Let $f=\log u$. Then we have
\begin{eqnarray}
\Delta f=-|\nabla f|^2-cu^{\alpha-1}.
\end{eqnarray}
Let $F=|\nabla f|^2+cu^{\alpha-1}$. Then we have $\Delta f=-F$ and by the well-known Weitzenbock-Bochner formula
$$\Delta|\nabla f|^2=2\nabla f \nabla\Delta f+2|\nabla^2f|^2+2{\rm Ric}(\nabla f, \nabla f),$$
where $\nabla^2f$ is the Hessian of $f$. Since $c>0$ and $\alpha\leq 1$,
we obtain by the above two inequalities
\begin{eqnarray*}
\Delta F&=&\Delta|\nabla f|^2+c\Delta u^{\alpha-1}
\\&=&-2\nabla f\nabla F+2|\nabla^2f|^2+2{\rm Ric}(\nabla f, \nabla f)
\\&+&c(1-\alpha)u^{\alpha-1}F+c(1-\alpha)^2u^{\alpha-1}|\nabla f|^2
\\&\geq&-2\nabla f\nabla F+\frac{2}{n}F^2-2KF
\end{eqnarray*}
on $B_p(2R)$, where we used the fact that
$|\nabla^2f|^2\geq\frac{1}{n}(\Delta f)^2$. Then following Yang's
proof line by line we finish the proof of theorem \ref{1-thm-0}.

\bibliographystyle{Plain}

\begin{thebibliography}{10}


\bibitem{Brighton2013}
K. Brighton, A Liouville-type theorem for smooth metric measure
spaces, J. Geom. Anal. 23 (2013), 562-570.



\bibitem{GidasSpruck-1981}
B. Gidas, J. Spruck, Global and local behavior of positive solutions of nonlinear elliptic equations, Comm. Pure Appl. Math. 34(1981), no.4, 525-598.


\bibitem{Guowei-2006}
Z. Guo, J. Wei, Hausdoff dimension of ruptures for solutions of a semilinear equation with singular
nonlinearity, Manuscripta Math. 120 (2006), 193-209.



\bibitem{HM-2017}
G. Y. Huang, B. Q. Ma, Gradient estimates and Liouville type theorems for a nonlinear
elliptic equation, Arch. Math. 105 (2015), 491-499.

\bibitem{HM-GeomDedicata2017}
G. Y. Huang, B. Q. Ma, Hamilton's gradient estimates of porous medium and fast
diffusion equations, Geom. Dedicata, 188 (2017), 1-16.

\bibitem{HLZ-2018}
G. Y. Huang, Z. Li, Liouville type theorems of a nonlinear elliptic equation
for the $V$-Laplacian, Anal. Math. Phys.  DOI 10.1007/s13324-017-0168-6.

\bibitem{Li-1991}
J. Y. Li, Gradient estimate for the heat kernel of a complete Riemannian manifold and its applications,
J. Funct. Anal. 97 (1991) 293-310.

\bibitem{Qian-2017}
B. Qian, Yau's gradient estimates for a nonlinear elliptic equation, Arch. Math. 108 (2017), 427-435.

\bibitem{schoenyau}
R. Schoen, S.-T. Yau, Lectures on Differential Geometry, International Press, 1994.

\bibitem{Yang2010}
Y. Y. Yang, Gradient estimates for the equation $\Delta u+cu^{-\alpha}=0$ on Riemannian manifolds, Acta
Math. Sin. (Engl Ser), 26 (2010), 1177-1182.

\bibitem{Yau75}
S. T. Yau, Harmonic functions on complete Riemannian manifolds, Comm. Pure Appl. Math. 28 (1975), 201-228.


\end{thebibliography}

\end{document}